\documentstyle[12pt]{article}
\begin{document}
\renewcommand{\thefootnote}{\fnsymbol{footnote}}
\newpage
\pagestyle{empty}
\setcounter{page}{0}
\renewcommand{\theequation}{\arabic{equation}}
\setcounter{equation}{0}
\newfont{\twelvemsb}{msbm10 scaled\magstep1}
\newfont{\eightmsb}{msbm8}
\newfont{\sixmsb}{msbm6}
\newfam\msbfam
\textfont\msbfam=\twelvemsb
\scriptfont\msbfam=\eightmsb
\scriptscriptfont\msbfam=\sixmsb
\catcode`\@=11
\def\Bbb{\ifmmode\let\next\Bbb@\else
  \def\next{\errmessage{Use \string\Bbb\space only in math mode}}\fi\next}
\def\Bbb@#1{{\Bbb@@{#1}}}
\def\Bbb@@#1{\fam\msbfam#1}
\newfont{\twelvegoth}{eufm10 scaled\magstep1}
\newfont{\tengoth}{eufm10}
\newfont{\eightgoth}{eufm8}
\newfont{\sixgoth}{eufm6}

\newfam\gothfam
\textfont\gothfam=\twelvegoth
\scriptfont\gothfam=\eightgoth
\scriptscriptfont\gothfam=\sixgoth
\def\frak{\frak@}
\def\frak@#1{{\fam\gothfam{{#1}}}}
\def\frak@@#1{\fam\gothfam#1}
\catcode`@=12
%
%
%
\def\CC{{\Bbb C}}
\def\NN{{\Bbb N}}
\def\QQ{{\Bbb Q}}
\def\RR{{\Bbb R}}
\def\ZZ{{\Bbb Z}}
\def\cA{{\cal A}}          \def\cB{{\cal B}}          \def\cC{{\cal C}}
\def\cD{{\cal D}}          \def\cE{{\cal E}}          \def\cF{{\cal F}}
\def\cG{{\cal G}}          \def\cH{{\cal H}}          \def\cI{{\cal I}}
\def\cJ{{\cal J}}          \def\cK{{\cal K}}          \def\cL{{\cal L}}
\def\cM{{\cal M}}          \def\cN{{\cal N}}          \def\cO{{\cal O}}
\def\cP{{\cal P}}          \def\cQ{{\cal Q}}          \def\cR{{\cal R}}
\def\cS{{\cal S}}          \def\cT{{\cal T}}          \def\cU{{\cal U}}
\def\cV{{\cal V}}          \def\cW{{\cal W}}          \def\cX{{\cal X}}
\def\cY{{\cal Y}}          \def\cZ{{\cal Z}}
\def\qed{\hfill \rule{5pt}{5pt}}
\def\id{\mbox{id}}
\def\ggo{{\frak g}_{\bar 0}}
\def\uqggo{\cU_q({\frak g}_{\bar 0})}
\def\uqggp{\cU_q({\frak g}_+)}
\def\typeA{{\em type $\cA$}}
\def\typeB{{\em type $\cB$}}
\newtheorem{lemma}{Lemma}
\newtheorem{prop}{Proposition}
\newtheorem{theo}{Theorem}
\newtheorem{Defi}{Definition}
$$
\;
$$
$$
\;
$$
$$
\;
$$
$$
\;
$$
\vfill
\vfill
\begin{center}

{\LARGE {\bf {\sf
The Gervais-Neveu-Felder equation for the Jordanian quasi-Hopf
$U_{h;y}(sl(2))$ algebra
}}} \\[2cm]

{\large  A. Chakrabarti$^{1}$\footnote{chakra@cpht.polytechnique.fr}
and R. Chakrabarti$^{2}$}\\
{\em
    $^{1}$ Centre de Physique Th\'eorique \footnote{Laboratoire
    Propre du CNRS UPR A.0014}, Ecole Polytechnique, 91128 Palaiseau
    Cedex, France\\
    $^{2}$ Department of Theoretical Physics, University of Madras,\\
    Guindy Campus, Madras 600 025, India\\
}
\end{center}
\smallskip

\smallskip

\smallskip

\smallskip

\smallskip

\vfill

\begin{abstract}
Using a contraction procedure, we construct a twist
operator that satisfies a shifted cocycle condition, and leads
to the Jordanian quasi-Hopf $U_{h;y}(sl(2))$ algebra. The
corresponding universal ${\cal R}_{h}(y)$ matrix obeys a
Gervais-Neveu-Felder equation associated with the
$U_{h;y}(sl(2))$ algebra. For a class of representations,
the dynamical Yang-Baxter equation may be expressed as a
compatibility condition for the algebra of the Lax operators.
\end{abstract}

\vfill

\vfill

$$
\;
$$
$$
\;
$$

\newpage

\pagestyle{plain}

\noindent
Recently a class of invertible maps between the classical $sl(2)$
and the non-standard Jordanian $U_h(sl(2))$ algebras  has been
obtained ~\cite{ACC96}-\cite{ACCS99}.  The classical
and the Jordanian coalgebraic
structures may be related ~\cite{ACC98}-\cite{KLM98}
by the twist operators
corresponding to these maps. Following the first twist leading from
the classical to the Jordanian Hopf structure, it is possible
to envisage a second twist leading to a quasi-Hopf quantization
of the Jordanian $U_h(sl(2))$ algebra. By explicitly constructing
the appropriate universal twist operator that satisfies a shifted
cocycle condition, we here obtain the Gervais-Neveu-Felder (GNF)
equation satisfied by the universal ${\cal R}$  matrix of a
one-parametric quasi-Hopf deformation of the
$U_h(sl(2))$ algebra.

The GNF equation corresponding to the standard Drinfeld-Jimbo
deformed $U_q(sl(2))$ algebra was studied in the context of
Liouville field theory ~\cite{GN84}, quantization of
Kniznik-Zamolodchikov-Bernard equation ~\cite{Fe94} and the
quantization of the Calogero-Moser model in the $R$ matrix
formalism ~\cite{BBB96}. The general
construction of the twist operators
leading to the GNF equation corresponding to the quasi-triangular
standard Drinfeld-Jimbo deformed $U_{q}({\sf g})$ algebras
and superalgebras were obtained in ~\cite{Fr97}-\cite{ABRR97}.

For the sake of completeness, we start by enlisting the general
properties of a quasi-Hopf algebra ${\cal A}$ ~\cite{K95}. For all
$a\in{\cal A}$ there exist an invertible element
${\Phi}\in {\cal A}\otimes {\cal A}\otimes {\cal A}$ and the
elements $(\alpha, \beta)\in {\cal A}$, such that
\begin{eqnarray}
(\hbox{id}\otimes\bigtriangleup)\bigtriangleup(a)
&=&\Phi\,\,(\bigtriangleup
\otimes \hbox{id})(\bigtriangleup(a))\,\,\Phi^{-1},\nonumber\\
(\hbox{id} \otimes \hbox{id} \otimes
\bigtriangleup)(\Phi)\,\,(\bigtriangleup
\otimes \hbox{id} \otimes \hbox{id})(\Phi)
&=&(1 \otimes \Phi)\,\,(\hbox{id}  \otimes \bigtriangleup
\otimes \hbox{id})(\Phi)\,\,(\Phi \otimes 1),\nonumber\\
(\varepsilon \otimes \hbox{id})\circ
\bigtriangleup &=& \hbox{id},\nonumber\\
(\hbox{id} \otimes \varepsilon)\circ
\bigtriangleup  &=& \hbox{id},\nonumber\\
\sum_{r} S(a_{r}^{(1)})\,\, \alpha \,\, a_{r}^{(2)} &=&
\varepsilon(a)\,\, \alpha,\nonumber\\
\sum_{r} a_{r}^{(1)}\,\, \beta\,\,S(a_{r}^{(2)})&=&
\varepsilon(a)\,\,\beta,\nonumber\\
\sum_{r}X_{r}^{(1)}\,\,\beta\,\,S(X_{r}^{(2)})\,\,\alpha\,\,
X_{r}^{(3)}&=&1,\nonumber\\
\sum_{r} S(\bar {X}_{r}^{(1)})\,\,\alpha \,\,
\bar {X}_{r}^{(2)}\,\,
\beta\,\,S(\bar {X}_{r}^{(3)})&=&1,
\label{eq:qhfdef}
\end{eqnarray}
where
\begin{equation}
\bigtriangleup(a)=\sum_{r} a_{r}^{(1)}\otimes a_{r}^{(2)},\quad
\Phi=\sum_{r} X_{r}^{(1)} \otimes X_{r}^{(2)} \otimes
X_{r}^{(3)},\quad
\Phi^{-1}=\sum_{r} \bar{X}_{r}^{(1)} \otimes  \bar{X}_{r}^{(2)} \otimes
\bar{X}_{r}^{(3)}.
\label{eq:apexp}
\end{equation}
A quasi-triangular quasi-Hopf algebra is equipped with a universal
${\cal R}$ matrix satisfying
\begin{eqnarray}
\bigtriangleup ^{op}(a)&=&{\cal R}\,\,\bigtriangleup(a)\,\,
{\cal R}^{-1},\nonumber\\
(\hbox{id} \otimes \bigtriangleup)({\cal R})&=&\Phi_{231}^{-1}\,\,
{\cal R}_{13}\,\,\Phi_{213}\,\,{\cal R}_{12}\,\,
\Phi_{123}^{-1},\nonumber\\
(\bigtriangleup \otimes \hbox{id})(\cal R)
&=&\Phi_{312}\,\,{\cal R}_{13}\,\,
\Phi_{132}^{-1}\,\,{\cal R}_{23}\,\,\Phi_{123}.
\label{eq:quatri}
\end{eqnarray}
The algebra is known as triangular if the additional relation
\begin{equation}
{\cal R}_{21}={\cal R}^{-1}
\label{eq:tri}
\end {equation}
is satisfied. In a quasi-triangular quasi-Hopf
algebra, the universal
${\cal R}$ matrix satisfies quasi-Yang-Baxter equation
\begin{equation}
{\cal R}_{12}\,\,\Phi_{312}\,\,{\cal R}_{13}\,\,\Phi_{132}^{-1}\,\,
{\cal R}_{23}\,\,\Phi_{123}=\Phi_{321}\,\,{\cal R}_{23}\,\,
\Phi_{231}^{-1}\,\,{\cal R}_{13}\,\,\Phi_{213}\,\,{\cal R}_{12}.
\label{eq:quayb}
\end{equation}
An invertible twist operator ${\cal F}\in{\cal A}\otimes{\cal A}$
satisfying the relation
\begin{equation}
(\varepsilon \otimes \hbox{id})({\cal F})=1=(\hbox{id} \otimes
\varepsilon)({\cal F})
\label{eq:fcoun}
\end{equation}
performs a gauge transformation as follows:
\begin{eqnarray}
\bigtriangleup_{\cal F}(a)&=&{\cal F}\,\bigtriangleup (a)\:
{\cal F}^{-1},\nonumber\\
\Phi_{\cal F}&=&{\cal F}_{23}(\hbox{id}
\otimes \bigtriangleup)({\cal F})
\,\,\Phi\,\,(\bigtriangleup \otimes \hbox{id})({\cal F}^{-1})
{\cal F}_{12}^{-1},\nonumber\\
\alpha_{\cal F}&=&\sum_{r}S({\bar f}_{r}^{(1)})\,\,
\alpha\,\,{\bar f}_{r}^{(2)},\nonumber\\
\beta_{\cal F}&=&\sum_{r} f_{r}^{(1)}\,\,\beta\,\,
S( f_{r}^{(2)}),\nonumber\\
{\cal R}_{\cal F}&=&{\cal F}_{21}{\cal R}{\cal F}^{-1},
\label{eq:ggtr}
\end{eqnarray}
where
\begin{equation}
{\cal F}=\sum_{r}f_{r}^{(1)}\otimes f_{r}^{(2)},\qquad
{\cal F}^{-1}=\sum_{r}{\bar f}_{r}^{(1)}\otimes{\bar f}_{r}^{(2)}.
\label{eq:fexp}
\end{equation}

The Jordanian Hopf algebra $U_{h}(sl(2))$ is generated by the
elements $(T^{\pm 1}\,(=e^{\pm hX}), Y, H)$, satisfying the
algebraic relations ~\cite{O92}
\begin{equation}
[H,T^{\pm 1}]=T^{\pm 2}-1,\,\,
[H,Y]=\,-\,\frac{1}{2}\,\left(Y(T+T^{-1})
\,+\,(T+T^{-1})Y \right),\,\,[X,Y]=H,
\label{eq:jalg}
\end{equation}
whereas the coalgebraic properties are given by ~\cite{O92}
\begin{eqnarray}
\bigtriangleup(T^{\pm 1})&=&T^{\pm 1}\otimes T^{\pm 1},\,\,
\bigtriangleup(Y)=Y\otimes T+T^{-1}\otimes Y,\,\,
\bigtriangleup(H)=H\otimes T+T^{-1}\otimes H,\nonumber\\
\varepsilon(T^{\pm 1})&=&1,\qquad\varepsilon(Y)=
\varepsilon (H)=0,\nonumber\\
S(T^{\pm 1})&=&T^{\mp 1},\quad S(Y)=-TYT^{-1},
\quad S(H)=-THT^{-1}.
\label{eq:jcoalg}
\end{eqnarray}
The universal ${\cal R}_{h}$ matrix of the triangular Hopf algebra
$U_{h}(sl(2))$ is given in a convenient form ~\cite{BH96} by
\begin{equation}
{\cal R}_{h}=\,\hbox{exp}(-hX\otimes TH)
\,\,\hbox{exp}(hTH\otimes X).
\label{eq:junvr}
\end{equation}
An invertible nonlinear map of the generating elements of the
$U_{h}(sl(2))$ algebra on the elements of the classical $U(sl(2))$
algebra plays a pivotal role in the present work.
The map reads ~\cite{ACC98}
\begin{equation}
T={\tilde T},\quad Y=J_{-}-\frac{1}{4}\,h^{2}
\,J_{+}\,(J_{0}^{2}-1),\quad
H=\,{(1+{(hJ_{+})}^{2}\,)}^{1/2}\,J_{0},
\label{eq:map}
\end{equation}
where ${\tilde T}=hJ_{+}\,+\,{(1+{(hJ_{+})}^{2})}^{1/2}$.
The elements $(J_{\pm},J_{0})$ are the generators of
the classical $sl(2)$ algebra
\begin{equation}
[J_{0},J_{\pm}]=\,\pm \,2\,J_{\pm},\qquad [J_{+},J_{-}]=J_{0}.
\label{eq:clalg}
\end{equation}
The twist operator specific to the map
~(\ref{eq:map}), transforming the trivial classical
$U(sl(2))$ coproduct structure to the non-cocommuting coproduct
properties ~(\ref{eq:jcoalg}) of the Jordanian
$U_{h}(sl(2))$ algebra, has been obtained ~\cite{ACCS99}, \cite{CQ99}
as a series expansion in powers of $h$. The
transforming operator between the two above-mentioned antipode
maps has been obtained ~\cite{CQ99} in a closed form.

Our present derivation of the GNF equation corresponding to the
Jordanian $U_{h}(sl(2))$ algebra closely parallels the description
in ~\cite{BBB96}. These authors obtained the solutions of the GNF
equation in the case of the standard Drinfeld-Jimbo deformed
quasi-Hopf $U_{q;x}(sl(2))$ algebra by constructing the
universal twist operator depending on a parameter $x$ :
\begin{eqnarray}
{\cal F}(x)&=&\sum_{k=0}^{\infty}(-1)^{k}\,\,
\frac{(q-q^{-1})^{k}}{[k]_{q}\,!}\,\,x^{2k}\,q^{k(k+1)/2}
\biggl[\,\, \prod_{l=1}^{k}{\left(1\otimes1\,-\,x^{2}q^{2l}
\,\,1\otimes q^{2{\cal J}_{0}}\right)}^{-1}
\biggr]\,\,\times\nonumber\\
& &\times\quad q^{\frac{k}{2}\,{\cal J}_{0}}\,\,
{\cal J}_{+}^{k}\,\,\otimes\,\,q^{\frac{3k}{2}
\,{\cal J}_{0}}\,\,{\cal J}_{-}^{k}\,,
\label{eq:qf}
\end{eqnarray}
where\,\,$[n]_{q}=\bigl(q^{n}-q^{-n}\bigr)/\bigl(q-q^{-1}\bigr)$.
The generators of the $U_{q}(sl(2))$ algebra
satisfies  ~\cite{K95} the relations
\begin{equation}
q^{{\cal J}_{0}}\,\,{\cal J}_{\pm}\,\,q^{-{\cal J}_{0}}
=q^{\pm 2}\,\,{\cal J}_{\pm},\quad
[{\cal J}_{+},\,{\cal J}_{-}]={[{\cal J}_{0}]}_{q}.
\label{eq:qalg}
\end{equation}

A key ingredient in our method is  the contraction technique
developed in ~\cite{ACC98}, where a matrix $G$
\begin{equation}
G=E_{q}(\eta {\cal J}_{+})\,\otimes \,E_{q}(\eta{\cal J}_{+}),
\quad \eta=\frac{h}{q-1}
\label{eq:conmat}
\end{equation}
performs a similarity transformation on the universal ${\cal R}_{q}$
matrix of the $U_{q}(sl(2))$ algebra ~\cite{K95}.
The twisted exponential $E_{q}({\chi})$ reads
\begin{equation}
E_{q}({\chi})=\sum_{n=0}^{\infty}\,
\frac{{\chi}^{n}}{[n]_{q}\,!}.
\label{eq:qexpo}
\end{equation}
The transforming matrix $G$ is singular in the $q\rightarrow 1$ limit.
The transformed $R_{h}^{j_{1};j_{2}}$ matrix for an arbitrary
$(j_{1};j_{2})$ represention
\begin{equation}
R_{h}^{j_{1};j_{2}}=\lim_{q \to 1} \left[G^{-1}\,
R_{q}^{j_{1};j_{2}}\,G \right]
\label{eq:rcont}
\end{equation}
is, however,  nonsingular and coincide, on account
of the map ~(\ref{eq:map}), with the result
obtained directly from the
expression ~(\ref{eq:junvr}) of the universal
${\cal R}_{h}$ matrix. In the
above contraction process the following two identities
play a crucial role:
\begin{eqnarray}
{(E(\eta\,{\cal J}_{+}))}^{-1}\,q^{\alpha\,{\cal J}_{0}\,/\,2}
\,E(\eta\,{\cal J}_{+})\,&=&\,{\cal T}_{(\alpha )}\,
q^{\alpha\,{\cal J}_{0}\,/\,2},\nonumber\\
{(E(\eta\,{\cal J}_{+}))}^{-1}\,{\cal J}_{-}\,
E(\eta\,{\cal J}_{+})&=& \,-\,{\frac {\eta}
{q-q^{-1}}}\,\bigl(\,{\cal T}_{(1)}\,
q^{{\cal J}_{0}}\,-\,{\cal T}_{(-1)}\,q^{-{\cal J}_{0}} \bigr)\,
+\,{\cal J}_{-},
\label{eq:ident}
\end{eqnarray}
where ${\cal T}_{(\alpha)}={(E(\eta\,{\cal J}_{+}))}^{-1}\,
E(q^{\alpha}\eta\,{\cal J}_{+})$.
In the $q\rightarrow 1$ limit, it may be proved ~\cite{ACC98}
\begin{equation}
\lim_{q\rightarrow 1}\,{\cal T}_{(\alpha)}={\tilde T}^{\alpha}
=T^{\alpha}.
\label{eq:tmap}
\end{equation}
The second equality in ~(\ref{eq:tmap}) follows from the
map ~(\ref{eq:map}).

Using the contraction scheme discussed above we now
obtain an one-parametric twist operator ${\cal F}_{h}(y)
\in U_{h}(sl(2))\,\otimes\,U_{h}(sl(2))$, which satisfies a shifted
cocycle condition. The twist operator ${\cal F}_{h}(y)$ gauge
transforms {\it{\`a} la} ~(\ref{eq:ggtr}) the Jordanian Hopf algebra

$U_{h}(sl(2))$ to a quasi-Hopf $U_{h;y}(sl(2))$ algebra and the
transformed universal ${\cal R}_{h}(y)$ matrix satisfies the
corresponding GNF equation. To this end we first compute
\begin{equation}
{\tilde {\cal F}}(y)=\lim_{q\rightarrow 1} \,\lgroup
G^{-1}\:{\cal F}(x) G\rgroup_{x^{2}=\,y\,(q-1)}\,,
\label{eq:fcontr}
\end{equation}
where ${\cal F}(x)$ is given by ~(\ref{eq:qf}).

A new feature here is the
reparametrization described by
\begin{equation}
y=\,{\frac {x^{2}}{q-1}},
\label{eq:repar}
\end{equation}
which is necessary for obtaining {\it nonsingular} result in the
$q\rightarrow 1$ limit. In ~(\ref{eq:repar}) we assume that
$x\rightarrow 0$ in the $q\rightarrow 1$ in such a way that $y$
remains finite. Following the above procedure in the
said limit we obtain
\begin{equation}
{\tilde {\cal F}}(y)= \sum_{k=0}^{\infty}\,
{\frac{{(hy)}^{k}}{k!}}\,\,
{\left({\tilde T}{J}_{+}\right)}^{k}\,\otimes\,
{\left({\tilde T}^{3}({\tilde T}-{\tilde T}^{-1})\right)}^{k}.
\label{eq:fnonsn}
\end{equation}
The rhs of ~(\ref{eq:fnonsn}) is interpreted
on account of the map ~(\ref{eq:map}) as
an element of $U_{h}(sl(2))\,\otimes\,U_{h}(sl(2))$. Identifying
this in the above sense with the twist operator ${\cal F}_{h}(y)
\left(\,={\tilde{\cal F}}(y)\right)$ we now
obtain the crucial result
\begin{equation}
{\cal F}_{h}(y)=\,\hbox{exp}\left({\frac{y}{2}}\,(1-T^{2})\,\otimes
\,(T^{2}-T^{4})\right).
\label{eq:hf}
\end{equation}
The above twist operator ${\cal F}_{h}(y)$ satisfies the property
~(\ref{eq:fcoun}). Following the arguments in ~\cite{BBB96}
we express ${\cal F}_{h}(y)$ as
a shifted coboundary
\begin{equation}
{\cal F}_{h}(y)=\bigtriangleup ({\cal M}(y))\,\,
\left(1\,\otimes \,{\cal M}^{-1}(y)\right)\,\,
\left({\cal M}^{-1}(y\,T_{(2)}^{4})\,\otimes\,1\right),
\label{eq:fshco}
\end{equation}
where the expression for the boundary reads
\begin{equation}
{\cal M}(y)=\,\hbox{exp}\left( {\frac {y}{2}}(1-T^{2})\right).
\label{eq:boun}
\end{equation}
The operator ${\cal F}_{h}(y)$ given by ~(\ref{eq:hf})
satisfies the following shifted cocycle condition
\begin{equation}
\left(1\,\otimes\,{\cal F}_{h}(y)\right)\,
\left[(\hbox{id}\,\otimes\,\bigtriangleup)\,{\cal F}_{h}(y)\right]=
\left({\cal F}_{h}\bigl( y\,T_{(3)}^{4}\bigr) \,\otimes\,1\right)\,
\left[(\bigtriangleup\,\otimes\,\hbox{id})\,{\cal F}_{h}(y)\right].
\label{eq:shcocl}
\end{equation}
Following ~(\ref{eq:ggtr}) the transformed
coproduct property may now be read as
\begin{equation}
\bigtriangleup_{y}(a)={\cal F}_{h}(y)\,\,\bigtriangleup (a)\,\,
{\cal F}^{-1}_{h}(y)\quad\hbox{for all}\,\,a\in U_{h;y}(sl(2)).
\label{eq:coptra}
\end{equation}
It may now be shown that the shifted cocycle condition is a
consequence of the following shifted coassociativity
property:
\begin{equation}
\left(\hbox{id}\,\otimes\,\bigtriangleup_{y} \right)\circ
\bigtriangleup_{y}(a)=\left( \bigtriangleup_{y\,T_{(3)}^{4}}
\,\otimes\,\hbox{id}\right) \circ\bigtriangleup_{y}(a).
\label{eq:shcoass}
\end{equation}

Following ~(\ref{eq:ggtr}) the gauge-transformed
universal ${\cal R}_{h}(y)$
matrix for the Jordanian quasi-Hopf $U_{h;y}(sl(2))$ algebra reads
\begin{equation}
{\cal R}_{h}(y)={\cal F}_{h\:\, 21}(y)\,{\cal R}_{h}\,
{\cal F}_{h}^{-1}(y).
\label{eq:rtran}
\end{equation}
The coassociator $\Phi (y)$ corresponding to the Jordanian
quasi-Hopf $U_{h;y}(sl(2))$ algebra may be obtained for the above
constuction of the twist operator obeying the shifted cocycle
condition ~(\ref{eq:shcocl}). Using
~(\ref{eq:ggtr}), (\ref{eq:hf}) and~(\ref{eq:shcocl}) we obtain
\begin{eqnarray}
\Phi (y)&=&{\cal F}_{h\:\, 12}(y\,T_{(3)}^{4})\:\:
{\cal F}_{h\:\, 12}^{-1}(y)\nonumber\\
&=&\hbox{exp}\;\left[\,-\,\frac{y}{2}\bigl(1-T^{2}\bigr)\,\otimes\,
\bigl(T^{2}-T^{4}\bigr)\,\otimes\,\bigl(1-T^{4}\bigr)\right].
\label{eq:coastr}
\end{eqnarray}
The elements $\alpha (y)$ and $\beta (y)$, characterizing
the antipode map of the $U_{h;y}(sl(2))$ algebra
may be similarly obtained from
~(\ref{eq:ggtr}), (\ref{eq:jcoalg}) and~(\ref{eq:hf}):
\begin{equation}
\alpha (y)=\,\hbox{exp}\,\left[\frac{y}{2}\,
{\bigl(1-T^{2}\bigr)}^{2}\right],\qquad
\beta (y)=\,\hbox{exp}\,\left[\,-\frac{y}{2}\,
{\bigl(1-T^{-2}\bigr)}^{2}\right].
\label{eq:albe}
\end{equation}
Using the guage transformation property of the universal ${\cal R}$
matrix in ~(\ref{eq:ggtr}) and our construction ~(\ref{eq:hf}) of
the twist operator, we now discuss the GNF equation associated with
the Jordanian quasi-Hopf $U_{h;y}(sl(2))$ algebra. The relations
~(\ref{eq:ggtr}),~(\ref{eq:hf}) and ~(\ref{eq:coastr}) lead to the
transformation property
\begin{equation}
{\cal R}_{h\:\,12}\,\left(y\, T_{(3)}^{4}\right)=\Phi_{213}(y)\,
{\cal R}_{h\:\,12}\,(y)\,\Phi_{123}^{-1}(y).
\label{eq:rshift}
\end{equation}
Now the quasitriangularity property of $U_{h;y}(sl(2))$ algebra
implies via ~(\ref{eq:quatri}),~(\ref{eq:coastr}) and
~(\ref{eq:rshift}) the following relations:
\begin{eqnarray}
\bigl(\hbox{id}\,\otimes\,\bigtriangleup_{y} \bigr)
\,{\cal R}_{h}\,(y)={\cal F}_{h\:\,23}(y)\,{\cal F}_{h\:\,23}^{-1}
\left(y\,T_{(1)}^{4}\right)\,{\cal R}_{h\:\,13}(y)\,
{\cal R}_{h\:\,12}\left(y\,T_{(3)}^{4}\right),\nonumber\\
\bigl(\bigtriangleup_{y}\,\otimes\,\hbox{id}\bigr)
\,{\cal R}_{h}\,(y)={\cal R}_{h\:\,13}\left(y\,T_{(2)}^{4}\right)\,
{\cal R}_{h\:\,23}(y)\, {\cal F}_{h\:\,12}
\left(y\,T_{(3)}^{4}\right)\,{\cal F}_{h\:\,12}^{-1}(y).
\label{eq:yquatr}
\end{eqnarray}
Using the transformation property~(\ref{eq:rshift}) we may now
recast the quasi Yang-Baxter equation ~(\ref{eq:quayb}) as the
GNF equation associated with the Jordanian quasi-Hopf
$U_{h;y}(sl(2))$ algebra:
\begin{equation}
{\cal R}_{h\:\,12}(y)\,
{\cal R}_{h\:\,13}\left(y\,T_{(2)}^{4}\right)\,
{\cal R}_{h\:\,23}(y)=
{\cal R}_{h\:\,23}\left(y\,T_{(1)}^{4}\right)\,
{\cal R}_{h\:\,13}(y)
{\cal R}_{h\:\,12}\left(y\,T_{(3)}^{4}\right).
\label{eq:hgnf}
\end{equation}

We now briefly consider the solutions of the above GNF equation
~(\ref{eq:hgnf}). Using the universal ${\cal R}_{h}(y)$ matrix
~(\ref{eq:rtran}), the twist operator ${\cal F}_{h}(y)$ in
~(\ref{eq:hf}) and the map ~(\ref{eq:map}) of the
generators of the $ U_{h}(sl(2))$ algebra
on the corresponding classical
elements, we may construct solutions of the GNF equation
~(\ref{eq:hgnf}). As illusrations we describe the
representions $R_{h}(y)$ for the ${\frac {1}{2}}\,\otimes\,j$
and the $1\,\otimes\,j$ cases. A $(2j+1)$ dimensional
representation of the classical $sl(2)$ algebra ~(\ref{eq:clalg})
\begin{eqnarray}
J_{+}\vert jm \rangle&=&(j-m)(j+m+1)\,\vert j\,m+1\rangle ,
\qquad J_{-}\vert jm \rangle=\vert j\,m-1\rangle,\nonumber\\
J_{0}\vert jm\rangle&=&m\, \vert jm \rangle,
\label{eq:clrep}
\end{eqnarray}
now, via the map ~(\ref{eq:map}), immediately furnishes the
corresponding $(2j+1)$ dimensional representation of the
$U_{h}(sl(2))$ algebra~(\ref{eq:jalg}). For the
$j={\frac {1}{2}}$ case, the generators remain undeformed.
For the $j=1$ case, we list the representation of
$U_{h}(sl(2))$ below.
\begin{eqnarray}
\lefteqn{(j=1)}\nonumber\\[0.2 cm]
& & X=\left(
      \begin{array}{ccc}
      0& 2& 0\\[0.2 cm]
      0& 0& 2\\[0.2 cm]
      0& 0& 0
      \end{array}\right),\qquad
      Y=\left(
      \begin{array}{ccc}
      0& {\scriptstyle {\frac{1}{2}}}h^{2}& 0\\[0.2 cm]
      1& 0& -{\scriptstyle {\frac{3}{2}}}h^{2}\\[0.2 cm]
      0& 1& 0
      \end{array}\right),\nonumber\\[0.2 cm]
& & H=\left(
      \begin{array}{ccc}
      2& 0& -4h^{2}\\[0.2 cm]
      0& 0& 0\\[0.2 cm]
      0& 0& -2
      \end{array}\right).
\end{eqnarray}
Using the above representations in the expression~(\ref{eq:rtran})
of the universal ${\cal R}_{h}(y)$ matrix,we obtain
\begin{equation}
R_{h}^{{\frac{1}{2}};\,j}\,(y)=\left(
\begin{array}{cc}
T& -hH+{\frac{1}{2}}\,h\,(T-T^{-1})
\Bigl(1+2y(1-T^{4})\Bigr)\\[0.2 cm]
0& T^{-1}
\end{array}\right)
\label{eq:rhfj}
\end{equation}
and
\begin{equation}
R_{h}^{1;\,j}\,(y)=\left(
\begin{array}{ccc}
T^{2}& A& B\\[0.2 cm]
0& 1& C\\[0.2 cm]
0& 0& T^{-2}
\end{array}\right),
\label{eq:ronj}
\end{equation}
where
\begin{eqnarray}
A&=&-2h\,TH\,-2hy\,\bigl(1-T^{2}\bigr)\,
\bigl(1-T^{4}\bigr),\nonumber\\
B&=&-2h^{2}\,\left[T^{2}-T^{-2}-2\,TH\bigl(1-T^{-2}\bigr)
-{(TH)}^{2}T^{-2}\right]\,
-4h^{2}y\,\bigl(1-T^{2}\bigr)\,\bigl(1+4T^{-2}-T^{4}\bigr)\nonumber\\
& &-4h^{2}y\, TH\,\bigl(1-T^{2}\bigr)\,\bigl(T^{2}-T^{-2}\bigr)
+2h^{2}y^{2}\,{\bigl(T-T^{-1}\bigr)}^{2}\,
{\bigl(1-T^{4}\bigr)}^{2},\nonumber\\
C&=&-2h\,\big(1-T^{-2}+THT^{-2}\bigr)
+2hy\,\bigl(1-T^{2}\bigr)\bigl(T^{2}-T^{-2}\bigr).
\label{eq:rabc}
\end{eqnarray}
{}From ~(\ref{eq:rhfj}) it follows that the
$R_{h}^{{\frac{1}{2}};{\frac{1}{2}}}$ matrix for the fundamental
$(1/2;1/2)$ case does not depend on the parameter $y$. The
$R_{h}(y)$ matrices for the higher representations, however,
nontrivially depend on $y$. The $R_{h}(y)$ matrices satisfy an
\lq\lq exchange symmetry" between the two sectors
of the tensor product spaces:
\begin{equation}
\left(R_{h}^{j_{1};\,j_{2}}(y)\right)_{km,\,ln}=
\left(R_{- h}^{j_{2};\,j_{1}}(y)\right)_{mk,\,nl}.
\label{eq:exchange}
\end{equation}

In the remaining part of the present work we recast the Jordanian
GNF equation ~(\ref{eq:hgnf}) as a compatibility condition for the
algebra of $ L$ operators. Using a new parametrization
$y=\hbox{exp}\,(z)$, we perform a translation
\begin{equation}
{\cal R}_{h\:\,12}(z)\rightarrow
{\cal R}_{h\:\,12}(z-2h\,X_{(3)})
\label{eq:rsym}
\end {equation}
to express ~(\ref{eq:hgnf}) in a symmetric form
\begin{eqnarray}
&&{\cal R}_{h\:\,12}(z-2hX_{(3)})\,
{\cal R}_{h\:\,13}(z+2hX_{(2)})\,
{\cal R}_{h\:\,23}(z-2hX_{(1)})\nonumber\\
&&\qquad\qquad={\cal R}_{h\:\,23}(z+2hX_{(1)})\,
{\cal R}_{h\:\,13}(z-2hX_{(2)})\,
{\cal R}_{h\:\,12}(z+2hX_{(3)}).
\label{eq:hgnfsm}
\end{eqnarray}
This is equivalent to the Jordanian GNF equation ~(\ref{eq:hgnf})
for the class of representations $\varrho_{j_{1};\,j_{2}}$
satisfying the property
\begin{equation}
\varrho_{j_{1};\,j_{2}}\Bigl(\Bigl[\bigl(X_{(k)}+X_{(l)}\bigr)\,
\partial_{z},{\cal R}_{h\:\,kl}(z)\Bigr]\Bigr)=0.
\label{eq:comrep}
\end{equation}
Adopting the procedure in ~\cite{BBB96} we here use the
following construction of the Lax operator for the
$U_{h;y}(sl(2))$ algebra
\begin{equation}
L_{13}(z)=\hbox{exp}\left[-2h\,\bigl(2X_{(1)}+X_{(3)}
\bigr)\,\partial_{z}\right]\,{\cal R}_{h\:\,13}(z)
\:\hbox{exp}\left[2h\,X_{(3)}\,\partial_{z}\right],
\label{eq:ldef}
\end{equation}
where the subscript $3$ denotes the quantum space. For the
representations satisfying ~(\ref{eq:comrep})
the relation ~(\ref{eq:hgnfsm}) may be expressed in a Lax martix
form
\begin{equation}
R_{h\:\,12}^{j_{1};\,j_{2}}(z-2h\,X_{(3)})\,L_{13}(z)\,L_{23}(z)=
L_{23}(z)\,L_{13}(z)\,R_{h\:\,12}^{j_{1};\,j_{2}}(z+2h\,X_{(3)}).
\label{eq:rellax}
\end{equation}
As illustrations we note that the representations
$R_{h}^{\frac{1}{2}\,;\,1}(z)$, $R_{h}^{1\,;\,\frac{1}{2}}(z)$
and $R_{h}^{1\,;\,1}(z)$ obtained from ~(\ref{eq:rhfj})
and~(\ref{eq:ronj}) satisfy the requirement
~(\ref{eq:comrep}).

To summarize, here we have constructed the Jordanian quasi-Hopf
$U_{h;y}(sl(2))$ algebra by explicitly obtaining the
relevant twist operator via a contraction method. In the
contraction method used here we start with the standard
Drinfeld-Jimbo deformed quasi-Hopf $U_{q;x}(sl(2))$ algebra and
use a suitable similarity transformation followed by a
$q\rightarrow 1$ limiting process. An important point here
is that the reparametrization as obtained in
~(\ref{eq:repar}) is essential for obtaining a
{\em nonsingular} twist operator  for the  $U_{h;y}(sl(2))$ algebra
in the $q\rightarrow 1$ limit. Our
contraction method has an advantage in that it furnishes
the dynamical quantities for the Jordanian
quasi-Hopf $U_{h;y}(sl(2))$ algebra from the corresponding
quantities of the standard Drinfeld-Jimbo deformed
quasi-Hopf $U_{q;x}(sl(2))$ algebra. The  present twist operator
associated with the $U_{h;y}(sl(2))$ algebra
satisfies a shifted cocycle condition. The universal
${\cal R}_{h}(y)$ matrix satisfies the GNF equation associated
with the $U_{h;y}(sl(2))$ algebra. For a special class of
representations, the GNF equation may be recast as a
compatibility condition of the $L$ operators. As an extension of
the present work, a similar formalism may be developed to describe
a quasi-Hopf quantization of the coloured Jordanian deformed
$gl(2)$ algebra considered in~\cite{Q97}, \cite{P98}, \cite{CQ99}.
A similar construction of the twist operators associated with
the quasi-Hopf deformation of the Jordanian $sl_{h}(N)$
algebra may also be attempted following the discussion in
\cite{ACC98}.

\vskip 1cm

\noindent {\bf Acknowledgments:}

\smallskip

One of us (RC) wishes to thank A. J. Bracken for a
kind invitation to the University of Queensland, where part of this
work was done.

\bibliographystyle{amsplain}

\begin{thebibliography}{}


\bibitem {ACC96} B. Abdesselam, A. Chakrabarti and R. Chakrabarti,
Mod. Phys. Lett.{\bf A11}\,(1996)\,2883.

\bibitem {ACC98} B. Abdesselam, A. Chakrabarti and R. Chakrabarti,
Mod. Phys. Lett.{\bf A13}\,(1998)\,779.

\bibitem {ACCS99}B. Abdesselam, A.Chakrabarti, R.Chakrabarti
and J. Segar, Mod. Phys. Lett.{\bf A14}\,(1999)\, 765.

\bibitem {CQ99} R. Chakrabarti and C. Quesne, Int. J. Mod. Phys.
{\bf A14}\,(1999)\,2511.

\bibitem {KLM98} P. P. Kulish, V. D. Lyakhovsky and A. I. Mudrov,
{\em Extended Jordanian twists for Lie algebras},
{\bf math.QA}/9806014.

\bibitem {GN84} J. L. Gervais and A. Neveu, Nucl. Phys. {\bf 238}\,
(1984)\,125.

\bibitem {Fe94} G. Felder, {\em Elliptic Quantum Groups},
Proc. ICMP, Paris (1994).

\bibitem {BBB96} O. Babelon, D. Bernard and E. Billey, Phys. Lett.
{\bf B375}\,(1996)\,89.

\bibitem {Fr97} C. Fronsdal, Lett. Math. Phys.
{\bf 40}\,(1997)\,117.

\bibitem {JKOS97} M. Jimbo, H. Konno, S. Odake and J. Shiraishi,
{\em Quasi-Hopf twistors for elliptic quantum groups},
{\bf q-alg}/9712029.

\bibitem {ABRR97} D. Arnaudon, E. Buffenoir, E. Ragoucy and
Ph. Roche, {\em Universal solutions of Quantum Dynamical
Yang-Baxter equations}, {\bf q-alg}/9712037.

\bibitem {K95} C. Kassel, {\em Quantum groups},\,(1995)\,
Springer Verlag.

\bibitem {O92} Ch. Ohn, Lett. Math. Phys. {\bf25}\,(1992)\,85.

\bibitem {BH96} A. Ballesteros and F. J. Herranz,
J. Phys. A: Math. Gen. {\bf 29}\, (1996)\, L311.

\bibitem {Q97} C. Quesne, J. Math. Phys. {\bf 38} (1997) 6018.

\bibitem {P98} P. Parashar, Lett. Math. Phys. {\bf 45} (1998) 105.

\end{thebibliography}

\end{document}